\documentclass{article}[10pt]

\usepackage{latexsym}
\usepackage{amssymb}
\usepackage{amsmath,amsfonts}
\usepackage{tikz}
\usepackage{graphicx}
\usepackage{multicol}
\usepackage{multirow}

\def \dd{\hfill \baseb \vskip .5cm}
\def \d{{\noindent \it Proof. } }

\def\MnB{{\cal M}_{n}(\B)}

\def \N {{\mathbb N}}

\def \R {{\mathbb R}}

\def \M {{\mathcal PSC}}

\def \0 {{\mathbf 0}}
\def \B {{\mathcal B}}
\def \al {{\alpha}}
\def \be {{\beta}}
\def \ga {{\gamma}}

\def \de {{\delta}}

\def \al {{\alpha}}
\newtheorem{theorem}{Theorem}[section]
\newtheorem{lemma}[theorem]{Lemma}
\newtheorem{cor}[theorem]{Corollary}
\newtheorem{defn}{Definition}[section]
\newtheorem{example}[theorem]{Example} %[section]

\newtheorem{pro}[theorem]{Proposition}
\newtheorem{rem}[theorem]{Remark}
\newtheorem{note}[theorem]{Note}

\def \be{{\mathbf e}}

\def\ex{\begin{example}}
\def\eex{\end{example}}
\def\exx{\end{example}}
\def\t{\begin{theorem}}
\def\tt{\end{theorem}}
\def\D{\begin{defn}}
\def\DD{\end{defn}}
\def\l{\begin{lemma}}
\def\ll{\end{lemma}}
\def\c{\begin{corollary}}
\def\cc{\end{corollary}}
\def\cj{\begin{conjecture}}
\def\cjj{\end{conjecture}}
\def\e{\begin{equation}}
\def\ee{\end{equation}}
\def\p{\begin{prop}}
\def\pp{\end{prop}}
\def\q{\begin{question}}
\def\qq{\end{question}}
\def\[{[\hskip-1pt [}
\def\]{]\hskip-1pt ]}

\def\B{{\mathbb B}}

\def\R{{\mathbb R}}

\def\N{{\mathbb N}}

\def\++{\boxplus}

\newcommand{\baseb}{\hfill \rule{2mm}{2mm}}

\begin{document}

\baselineskip 15pt

\title{Exponents of Primitive Symmetric Companion Matrices}
\date{}
\author {Monimala Nej and A. Satyanarayana Reddy\\
Department of 
Mathematics \\ Shiv Nadar 
University, India-201314 \\ (e-mail: 
mn636@snu.edu.in, satyanarayana.reddy@snu.edu.in).
  }
\maketitle
\begin{abstract}
A {\it symmetric companion matrix} is a matrix of the form $A +A^T$ where $A$ is a
companion matrix all of whose entries are in $\{0,1\}$ and $A^T$ is the transpose of $A.$ In this paper, we find the total number of primitive and the total number of imprimitive symmetric companion matrices. We establish formulas to compute the exponent of every primitive symmetric companion matrix. Hence the exponent set for the class of primitive symmetric companion matrices is completely characterized. We also obtain the number of primitive symmetric companion matrices with a given exponent for certain cases.
\end{abstract}
{\bf{AMS(2010)}}: 05C50, 05C38, 15B99.
%%%%%%%%%%%%%%%%%%%%%%%%%%%%%%%%
\section{Introduction}  
Let $n$ be a positive integer, and let $A, B \in M_n(\R),$ the set of all $n \times n$ real
matrices. We denote $ij$-th entry of $A$ by $a_{i,j}.$  If all of the entries of $A$ are positive, then $A$ is called a {\it positive matrix}, which is denoted by $A > 0.$ If all of the entries of $A$ are nonnegative,
then $A$ is called a {\it nonnegative matrix}, which is denoted by $A \geq 0.$ By $A > B$
(resp., $A \geq B$), we mean $A-B > 0$ (resp., $A-B
 \geq 0);$ that is, $a_{i,j} > b_{i,j}$
(resp., $a_{i,j} \geq b_{i,j}$) for all $i$ and $j.$
A matrix $A$ is called {\it primitive} if it is nonnegative and there exists a positive integer $k$ such that $A^k > 0.$ For a primitive matrix $A,$
the smallest such positive integer $k$ is called the {\it primitive exponent} of $A,$ or
simply, the {\it exponent} of $A,$ and is denoted by $exp(A).$ An irreducible
matrix that is not primitive is called {\it imprimitive}. For additional information
on irreducible matrices, see H. Minc \cite{minc}.

Let $\B$ denote the binary Boolean semiring, that is $\B$ is the set $\{0,1\}$ with arithmetic the same as for the reals, except that $1+1=1$.  The set $\MnB$ is the set of $n\times n$ matrices with entries in $\B$ which forms a semimodule with the usual definitions of addition and multiplication.

Let $A$ be a nonnegative matrix.  The {\em support} of $A$ (sometimes called the nonzero pattern of $A$), denoted $\overline{A}$, is the matrix in $\MnB$ such that $\overline{a_{i,j}}=0$ if $a_{i,j}=0$, and   $\overline{a_{i,j}}=1$ if $a_{i,j}$ is nonzero. Since the product or sum of two positive numbers is always positive and the magnitude of the nonzero entries in a primitive matrix is not important, we have that a nonnegative matrix $A$ is primitive if and only if $\overline{A}$ is primitive and  $exp(A) = exp(\overline{A})$.

Let $X$ be a nonempty subset of $M_n(\R).$ The {\it exponent set} of $X,$ denoted
$E(X),$ is given by
$$E(X) = \{k\in \N : \mbox{ there exists a primitive $A\in X$ with 
$exp(A)=k$}\}.$$ 
Note that $X$ can contain both primitive and imprimitive matrices, and also that
$E(X)$ is empty when $X$ contains no primitive matrices.

For positive integers $a$ and $b$ with $a \leq b,$ the set of all positive integers $k$
with $a \leq k \leq b,$ will be denoted by $\[a, b\].$ This set partitions into two sets, the
subset of all even integers in $\[a, b\],$ denoted by $\[a, b\]^e ;$ and the subset of all odd
integers in $\[a, b\],$ denoted by $\[a, b\]^o.$

Every positive matrix is  a primitive matrix with exponent one, that is, $1\in 
E(M_n(\R)).$
In 1950, Helmut Wielandt \cite{Wie} proved that if $A\in M_n(\R)$ is primitive, then  
$exp(A)\leq \omega_{n}$, where 
$\omega_{n}=(n-1)^2+1$. This bound is known as the {\bf Wielandt bound}.
For the proof, see Hans Schneider \cite{HS},  Holladay and 
Varga \cite{Ho:Va}, or  Perkins \cite{Per}. 
%or in  the book by R.A. Brualdi and Herbert J. Ryser \cite{bru}). 
With this remarkable result we have $E(M_n(\R))\subseteq \[1,\omega_n\].$ In 1964,
Dulmage and Mendelsohn \cite{D:M} showed that $E(M_n(\R))\subsetneq 
\[1,\omega_n\].$

%is attained only by the Wielandt graph in Figure \ref{fig:figure1}. The adjacency matrix is a primitive companion matrix (see below).

 The study of $exp(A)$ has been focusing on the following problems.
\begin{description}
\item[(MEP)]
{\it The maximum exponent problem}, i.e., to estimate the upper bound of $E(X)$ for a particular class $X$ of primitive matrices. Here people generally looks for the `best possible' upper bound, i.e., upper bound $k$ of $E(X)$ such that $k \in E(X).$ Sometimes it is called the exact upper bound. 
\item[(SEP)]
{\it The set of exponents problem}, i.e., to determine $E(X).$
\item[(EMP)]
{\it The extremal matrix problem}, i.e., to determine primitive matrices in $X$ with the maximum exponent.
\item[(SMP)]
The extremal matrix problem can be generalized further as {\it the set of matrices problem}, i.e., for any $k \in E(X),$ determine primitive matrices in $X$ with the exponent $k.$ 
\end{description} 
  
\begin{figure}[htb]
\begin{center}
\includegraphics[width=.5\linewidth, width=.5\textwidth]{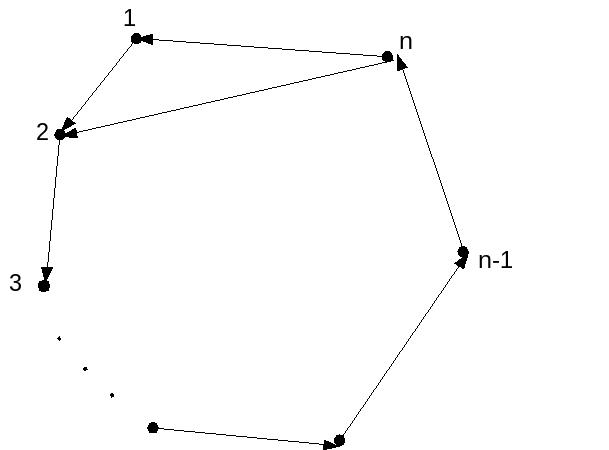}
\end{center}
\caption{The Wielandt graph $W_n.$}  \label{fig:figure1}
\end{figure}

Significant literature is available for these types of problems except for {\bf (SMP)}, the set of matrices problem. The papers \cite{Bru:Ross}, \cite{Lewin}, \cite{M:Y}, \cite{Liu1},\cite{B:M:W},\cite{Ross}, \cite{ZKM} studied  exponent sets for different classes of primitive matrices. In 1987, Shao \cite{shao} investigated the exponent set of $S_n,$ the class of $(0,1)$ primitive symmetric matrices. He showed that $E(S_n)=\[1,2n-2\] \setminus \[n,2n-1\]^o.$ Later, Liu et al. \cite{B:M:W} established the results for the primitive matrices in $S_n$ with zero trace. Here we  consider another class of primitive symmetric matrices of order $n.$ We denote this class by $\M_n,$ the class of primitive symmetric companion matrices. We show that the exponent set for this class is $\[2,2n-2\]^e=\[1,2n-2\] \setminus \[1,2n-1\]^o.$  Theorem $4$  of Fuyi et al. \cite{F:M:J} and the Wielandt graph in Figure \ref{fig:figure1} motivate us for such a consideration. Surprisingly, upper bounds as well as even numbers in the exponent sets given in \cite{shao} and \cite{B:M:W} are attained by  matrices that  belong to $\M_n.$ Not only that, in this paper we find the exponent of each $A \in \M_n.$ Hence we completely solve the set of matrices problem {\bf (SMP)} for $\M_n.$ Characterizations of matrices in $S_n$ with exponent $x \geq n-3$ can be found in literature. For instance, see Lichao and Cai \cite{Lic:Cai}. Furthermore, we obtain the number of matrices which attain the exact upper and lower bound of the exponent set. We refer readers to Theorem  6, 7, and 8 in \cite{D:M} where the authors found the number of primitive matrices with a given exponent. In 2007, Liu et al. \cite{Liuu:Youu:Yuu} found the number of  $n \times n$ $(0,1)$ primitive matrices with exponent $(b-1)^2+1,$ where $b$ is the Boolean rank of the matrix. Research in the area of the set of matrices problem and related counting problem has not seen significant progress in  recent years. This problem appears to be harder. The study of this problem is just the beginning, see Kim et al. \cite{Kim:Song:Hwang}.

\section{Notation}
Let $C_n=\left\{A\in M_n(\R): A=\begin{bmatrix}
0&1&0&0& \dots& 0\\
0&0&1&0& \dots& 0\\
\vdots &\vdots&\vdots&\vdots&\vdots &\vdots\\
0&0&0&0& \dots& 1\\
a_{0}&a_{1}&a_{2}&a_{3}& \dots& a_{n-1}
\end{bmatrix}, a_i\in \{0,1\}\right\}.$

That is $C_n$ is the  set of all
companion matrices of polynomials  of the form 
$x^{n}-\sum\limits_{i=0}^{n-1}a_{i}x^{i}$, where $a_i\in \{0,1\}.$ Since the first 
$n-1$ rows of every matrix in $C_n$ are fixed, hence it is sufficient to 
specify 
the last row. Thus  
there is a bijection between $C_n$ and $B_n$, where $B_n$ denotes the set 
of 
all binary strings of length $n,$ in particular 
$|C_n|=|B_n|=2^n,$ where $|S|$ denotes the cardinality of the set $S.$
Suppose $F:M_n(\R)\to M_n(\R)$ is given by $F(A)=A+A^T,$ where $A^T$ is the transpose of 
$A.$ Then every entry of $F(A)$, where $A \in C_n$ is either $0$ or $1$ except that 
$f_{n, n-1}=f_{n-1, n} \in \{1,2\}$ and that $f_{n, n}\in \{0,2\},$  where $(f_{i,j})=F(A)$.

If $A\in C_n$, then $A$ is called a {\it $(0,1)$ companion matrix} and we call $F(A)$ a {\it  symmetric companion matrix}.  
For $\al,\be \in \{0,1\}$, we define 
$$ C_n^{\al,\be}=\{F(A) : A \in C_n,  a_{n1}=\al, 
a_{n n}=\be\},$$ $$\M_n^{\al,\be}=\{ B \in C_n^{\al,\be} : \mbox{$B$ is primitive} 
\}.$$ 
Elements in $C_n^{\al,\be}$ will be denoted by $A_Y$, where $Y \in B_{n-3}$ and 
last row of $A_Y$ will be the row vector $[\al, Y,1, \be]$ or $[\al, Y,2, \be]$. 
Hence corresponding to each $Y$ there are two elements in $C_n^{\al,\be}$ which 
are both imprimitive or primitive simultaneously 
and in case of primitivity, both of them have the same exponent. 
\begin{rem}\label{rem:nn-1}
The number of imprimitive symmetric companion matrices of order $n$ or the number of matrices in $\M_n$ with a given exponent is always even.
\end{rem}
\begin{example}\label{ex:all}
For $n=3$, $B_3=\{000,001,010,100,011,101,110,111\},$ and
 $$C_3=\left\{ \begin{bmatrix}
               0& 1& 0\\
               0 &0 &1\\
               0&0&0
              \end{bmatrix}, \begin{bmatrix}
               0& 1& 0\\
               0 &0 &1\\
               0&0&1
              \end{bmatrix}, \ldots, \begin{bmatrix}
               0& 1& 0\\
               0 &0 &1\\
               1&1&1
              \end{bmatrix}\right\}.$$
  Since  the matrices $\begin{bmatrix}
               0& 1& 0\\
               1 &0 &1\\
               0&1&0
              \end{bmatrix}, \begin{bmatrix}
               0& 1& 0\\
               1 &0&2\\
               0&2&0 
               \end{bmatrix}$ are imprimitive, then $\M_3^{0,0}=\emptyset.$
              
  Also, $\M_3^{0,1}= \left\{ \begin{bmatrix}
               0& 1& 0\\
               1 &0 &1\\
               0&1&2
              \end{bmatrix}, \begin{bmatrix}
               0& 1& 0\\
               1 &0&2\\
               0&2&2 
               \end{bmatrix}
              \right\},  \,\,
 \M_3^{1,0}=\\ \left\{ \begin{bmatrix}
               0& 1& 1\\
               1 &0 &1\\
               1&1&0
              \end{bmatrix}, \begin{bmatrix}
               0& 1& 1\\
               1 &0&2\\
               1&2&0 
               \end{bmatrix}
              \right\}, \, {\text {and}}\,\, \; 
   \M_3^{1,1}= \left\{ \begin{bmatrix}
               0& 1& 1\\
               1 &0 &1\\
               1&1&2
              \end{bmatrix}, \begin{bmatrix}
               0& 1& 1\\
               1 &0&2\\
               1&2&2 
               \end{bmatrix}
              \right\}.$\\
 Moreover, it is easy to check that $E(\M_3^{0,1})=\{4\}$ and       
 $E(\M_3^{1,0})=\{2\}=E(\M_3^{1,1}).$ 
\exx

 To the best of our knowledge no one has studied the number of primitive and imprimitive 
matrices of a given class and also the  number of matrices with a given 
exponent 
from that class. Furthermore, there is no specific formula for 
 computing the exponent of a given matrix from a given class.  We address here these 
problems for the class of primitive symmetric companion matrices.

Here is an overview of the content of the paper. Our interest is in finding 
\begin{enumerate} \label{probsta}
 \item \label{probsta:1} $|\M_n|$, where $\M_n=\M_n^{0,0}\cup \M_n^{0,1}\cup 
\M_n^{1,0}\cup \M_n^{1,1},$ 
 \item \label{probsta:2}  $E(F(C_n))$ by finding $E(\M_n^{\al,\be})$, for 
$\al,\be \in \{0,1\}$ as   $E(F(C_n))=E(\M_n^{0,0})\cup E(\M_n^{0,1})\cup 
E(\M_n^{1,0})\cup E(\M_n^{1,1}),$
 \item \label{probsta:3}   $N_n^{\al,\be}(b)$, where $N_n^{\al,\be}(b)=|\{A\in 
\M_n^{\al,\be} : exp(A)=b, b\in E(\M_n^{\al,\be})\}|.$
\end{enumerate}

In Section~\ref{sec:loop} we will address Problem~\ref{probsta:1} that is 
we will find $|\M_n|$. 
The set $E(\M_n^{0,1}\cup \M_n^{1,1})$  and the numbers $N_n^{0,1}(b), 
N^{1,1}(b)$ will  be presented in Section~\ref{sec:loop}.
In Sections~\ref{sec:CBNL} and \ref{sec:NCNL} we will evaluate 
$E(\M_n^{1,0}), N_n^{1,0}(b)$ and $E(\M_n^{0,0}), N_n^{0,0}(b)$ respectively 
for 
specific values of $b$. At the end of Section \ref{sec:NCNL} we display a table in which $N_n^{\al, \be}(b), \; \al, \be \in \{0,1\}$ can be found for a small $n$ and for all possible $b.$
In the rest of this section we will see a few preliminaries and 
notations required for rest of the paper.

For each nonnegative matrix $A\in M_n(\R)$ we associate a digraph 
$D(A)$ with vertex set $V(A)\;(\mbox{or simply $V$})=\[1,n\]$ 
and 
edge set $E(A)=\{(r,s)\in V\times V : a_{rs}>0\}.$ There are an infinite number of  
matrices associated with a single digraph.  Further if the matrix $A$ is 
symmetric, then $D(A)$ is an
undirected graph or simply a  graph. In particular,  $D(F(A))$ is a graph for 
every $A\in C_n.$ If $X$ is a digraph, then  $X$ is a  primitive (imprimitive) 
if its adjacency matrix is  primitive (imprimitive)  and $exp(X)$ is defined to 
be  the exponent of 
its adjacency matrix.

For $A\in C_n^{\al,\be}$ we define $V_1(A), V_2(A)$ ( or simply $V_1,V_2$ if 
$A$ is clear form the context)  by 
$$V_1=\{i \in \[1,n-1\] : a_{ni}=0\}\;\; and\;\; V_2=\{i \in \[1,n-1\] : 
a_{ni}>0\}.$$ 
For  $U\subseteq \[1,n\]$ write  $U=U_1\cup U_2\cup \dots \cup U_r$, where for 
each $ k\in \[1,r\]$, $U_k=\[i_k, j_k\]$ with $1 \leq i_k \leq j_k \leq n,$ and for each $k \in \[2,r\],$ $j_{k-1}+2 \leq i_k.$ We define 
$m(U)=\max\{|U_1|, 
|U_2|,\ldots, |U_r|\}.$  For $U= \emptyset$, $m(U)$ assume to be zero.  
For example if $U=\{2,3,5,7,8,9,10,13,14\},$  then $U=\{2,3\}\cup \{5\}\cup 
\{7,8,9,10\}\cup \{13,14\}$ and 
$m(U)=|\{7,8,9,10\}|=4.$ From the construction of $C_n^{\al,\be}$, we have  
$n-1\in V_2,$  and thus   $0\le m(V_1)\le n-2.$ 

\ex\label{ex:1011}
Consider the graph with $V=\[1, 10\]$ in Figure \ref{fig:figure2}. For this, we have  $V_1=\{2\}\cup \{4\}\cup \{6\}\cup \{8\}$, $m(V_1)=1$ and $V_2=\{1,3,5,7,9\}.$ 
\exx
\begin{figure}[htb]
\begin{center}
\includegraphics[scale=.7]{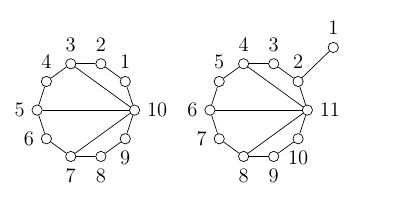}
\end{center}
\caption{Graphs of imprimitive symmetric companion matrices} \label{fig:figure2} 
\end{figure}

\section{Number of primitive symmetric companion matrices and 
$E(\M_n^{0,1}\cup \M_n^{1,1})$}\label{sec:loop}
The following theorem gives the number of primitive symmetric companion 
matrices. Before continuing  
recall that  a connected graph is a primitive graph
if and only if it has a cycle of odd length.  From the construction, $D(F(A))$ 
is a connected  graph for every $A\in C_n$ as  it 
 contains the path graph $1-2-3-\dots-n$ as a subgraph. Hence in order   to 
find $|\M_n^{\al,\be}|$ it is sufficient to 
 find  $|C_n^{\al,\be}\setminus \M_n^{\al,\be}|.$ That is we count the number 
of 
 matrices  $A$ in $C_n$ such that $a_{n1}=\al, a_{nn}=\be$ 
 and every cycle in  $D(F(A))$ is of even length. Further,  the two positions 
of 
 last row of  every matrix in  $C_n^{\al,\be}$
 are fixed, hence $|C_n^{\al,\be}|=2^{n-2}$ or equivalently 
$|\M_n^{\al,\be}|=2^{n-2}-|C_n^{\al,\be}\setminus\M_n^{\al,\be}|.$
 \pagebreak
\begin{theorem}\label{thm:NPSCM} Let $n\ge 4,$ and $\al\in \{0,1\}.$ Then
\begin{multicols}{2}
\begin{enumerate}
 \item\label{thm:NPSCM:1} $|{\M}_{n}^{\al,1}|=2^{n-2}.$ 
 \item\label{thm:NPSCM:2} $|{\M}_{n}^{1,0}|= \begin{cases}
                                                     2^{n-2} & \mbox{when $n$ 
is 
odd},\\
                                                     2^{n-2}-2^{\frac{n-2}{2}} 
& 
\mbox{when $n$ is even}.
                                                    \end{cases}$
  \item\label{thm:NPSCM:3} 
$|{\M}_{n}^{0,0}|=2^{n-2}-2^{\lfloor{\frac{n-1}{2}}\rfloor}$. 
\end{enumerate}
\end{multicols}
\end{theorem}
\d
Proof of Part \ref{thm:NPSCM:1}. Let $B=[b_{i,j}]\in C_n^{\al,1}.$ Then 
$B\in \M_{n}^{\al,1}$, as  there is a loop at the vertex $n$  in 
$D(B)$, which is  of odd length. Consequently  $C_n^{\al,1}\setminus 
\M_n^{\al,1}=\emptyset.$

Proof of Part \ref{thm:NPSCM:2}.
\begin{description}
 \item[$n$ is odd.]  Similar to proof of  Part \ref{thm:NPSCM:1}, now there is 
cycle of 
length $n$ in $D(B)$ for every $B\in C_n^{1,0}.$
 \item[$n$ is even.]  Let $A\in C_n^{1,0}$. Then $1,n-1\in V_2(A)$ as $n-1$ 
always belongs to $V_2(A)$. 
  If $V_2 \cap \[1,n\]^e\ne \emptyset$ then $D(A)$ contains an odd cycle hence 
$A\in \M_n^{1,0}.$ 
  Thus $A\in C_n^{1,0}\setminus \M_n^{1,0}$ if and only if $V_2\setminus 
\{1,n-1\} \subseteq \{3,5,\ldots,n-3\}.$
  There are $2^{\frac{n-4}{2}}$ such cases. Hence the result follows from 
Remark \ref{rem:nn-1}. 
  \end{description}  
 
 Proof of Part \ref{thm:NPSCM:3}. In this case,  $A\in C_n^{0,0}\setminus 
\M_n^{0,0}$ if and only if 
 $V_2\setminus \{n-1\}\subseteq \begin{cases} 
     \{3,5,7, \ldots ,n-3\}  & \mbox{when $n$ is even,}\\
     \{2,4,6, \ldots ,n-3\} &  \mbox{when $n$ is odd.}                        
                       \end{cases}$\\
The remainder of the proof is similar to the proof of the even case  of Part \ref{thm:NPSCM:2}.
\dd

\ex
 Both the graphs in Figure \ref{fig:figure2} are imprimitive.
\begin{enumerate}
 \item Consider the graph with $V=\[1,10\]$ in  Figure \ref{fig:figure2}. Then 
by 
using Remark \ref{rem:nn-1}
 it is easy to see that  
 the number of imprimitive symmetric companion matrices of 
order $10$ is $2^{(\frac{10-2}{2})} = 16$. These are 
the adjacency matrices of subgraphs obtained 
 by removing one or more  edges from $\{\{10,3\}, \{10,5\}, \{10,7\}\}$ of that 
graph. 
\item Similarly from the graph with $V=\[1, 11\]$ in Figure \ref{fig:figure2} it is easy to 
see 
that the number of imprimitive symmetric companion 
matrices of  order $11$ is  $2^{\lfloor\frac{11-1}{2}\rfloor} = 32.$
\end{enumerate}
\exx
%%%%%%%%%%%%%%%%%%%%%%%%%%%%%%%%%%%%%%%
%%%%%%%%%%%%%%%%%%%%%%%%%%%%%%%%%%%%%%%

Henceforth we restrict our investigation to  primitive symmetric companion  matrices, and 
focus on finding the exponent set of these matrices. In the paper \cite{shao} it is shown that the 
exponent of a primitive symmetric matrix is at most $2(n-1).$ Thus if $A$  is a primitive symmetric companion of order $n\ge 2,$  then $exp(A)\in \[2,2n-2\].$

 We now will give a  general procedure to find exponent of a primitive matrix.

Let $A$ be a primitive matrix and let  $exp(A:i,j)$ be the 
smallest positive integer $k$ such that there exists a walk of length $\ell$ from 
vertex $i$ to vertex $j$ in $D(A) \; \hbox{ for all }\; \ell \geq k. $   Equivalently, for $B=A^\ell$, 
$b_{i,j}>0 \hbox{ for all }  
\ell  \geq k.$
Let  $exp(A:i)$ is the smallest positive integer 
$p$ such that there exists a walk of length $p$ from vertex $i$ to any vertex 
$j$ in $D(A).$ Equivalently, every entry in the $i^{th}$ row of $A^p$ is 
positive. As a consequence, every entry in 
the $i^{th}$ row of $A^{p+1}$ is also  positive.

The following result can be found in Brualdi and Ryser \cite{bru}.
\begin{lemma}\label{lemma:exp(A)}
 Let $A\in M_n(\R)$ be a primitive matrix. Then $exp(A)=\max\limits_{1 \leq i,j \leq n}exp(A:i,j)=\max\limits_{1 \leq i \leq 
n}exp(A:i).$
\end{lemma}

\ex Consider the Fan graph, $F_8$, in  Figure \ref{fig:figure3}, whose adjacency matrix is a symmetric companion matrix.  It is easy to check  that   $exp(F_8:8,1)=1$,  
$exp(F_8:8,8)=2$.   
Also note that $exp(F_8:i)=2$ for $1 \leq i \leq n.$ Hence $exp(F_8)=2.$
\exx
\begin{figure}[htb]
\begin{center}
\includegraphics[]{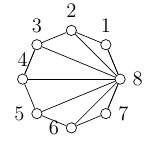}
\end{center}
\caption{The Fan graph $F_8$}  \label{fig:figure3}
\end{figure}

Hence our objective is to find  $exp(A:i,j)$ for all $i,j\in \[1,n\].$  In 
this paper, a  cycle at a vertex $i\in \[1,n\],$ we mean a walk from $i$ to $i,$ 
where 
repetition of internal vertices is allowed. 
Where as an elementary cycle or circuit we mean a cycle and  repetition of 
vertices is  allowed only at initial and terminal vertices. 
If there is a loop at a vertex $i$, then $exp(A:i,i)=1.$  
For the remaining cases, see Proposition $4$ in Liu et al.  \cite{B:M:W}. It shows that 
$exp(A:i,j)=d(i,j)+2k$ for some nonnegative integer $k$.  We include the proof for the sake of completeness.

\begin{pro}[\cite{B:M:W}]\label{pro:expoi}
 Let $A\in M_n(\R)$ be a symmetric  primitive matrix  and $i,j\in \[1,n\].$ Then 
$exp(A:i,j)=\ell-1$, 
 where $\ell$ is  the   length of the shortest walk  from
 $i$ to $j$ in $D(A)$ such that $2\nmid d(i,j)+\ell,$ and  $d(i,j)$ is the 
distance between  $i$ and $j$ in the graph $D(A).$ In particular, if 
$a_{ii}=0$, 
then  $exp(A:i,i)=\ell-1,$ where $\ell$ is the length of 
 shortest odd cycle from $i$ to $i.$
\end{pro}
\d
 If  $A\in M_n(\R)$ is a symmetric  primitive matrix,  then as discussed 
earlier 
$D(A)$ is a connected graph. 
 As a consequence  $ij^{th}$ entry of  $A^{d(i,j)+2k}$ is positive for every 
nonnegative integer $k.$ Suppose $\ell$ is the 
 length of the shortest walk from $i$ to $j$ with  parity different from that 
of $d(i,j)$, then $ij^{th}$ entry of  $A^{(\ell-1)+t}$ is positive for every 
nonnegative 
integer $t.$ Hence the result follows.
\dd

\begin{cor} [\cite{B:M:W}]\label{cor:expoi}
Let $G$ be a primitive simple graph, and let $i$ and $j$ be any pair of 
vertices in $V(G)$. 
If there are two walks $P_1, P_2$ from $i$ to $j$ with lengths $k_1$ and $k_2$ 
respectively, where $2\nmid k_1+k_2,$ 
then $exp(G:i,j) \leq max\{k_1,k_2\}-1.$ 
\end{cor}

The following theorem provides  a necessary condition for the exponent of primitive 
symmetric matrix  to be odd.

\begin{theorem}[\cite{F:M:J}]
 Let $X$ be a  primitive graph and suppose that $exp(X)$ is odd. Then $X$ contains two 
vertex  disjoint odd cycles.
\end{theorem}
From the construction of sets $\M_n^{\al,\be}$, if $A\in \M_n^{\al,\be}$, then 
every cycle in $D(A)$  contains  the vertex $n.$ 
Hence the following result.
\begin{cor}\label{cor:even}
 If $A\in \M_n^{\al,\be}$, then $exp(A)$ is even.
\end{cor} 

Thus, if $d(i,j)$ is odd, then $exp(A:i,j)<exp(A).$ Before proceeding to next result, we need 
the following notation.
\D
Let $B_n^{q,k}\subseteq B_n$ denote the 
set of all binary strings of length $n$ with $q$ zeros and having at least 
one  longest subword of zeros of length $k.$ \DD
For example, 
$B_6^{4,2}=\{100100,010100,010010,001100,001010,001001\}.$ 
Consequently, a  necessary condition for  $B_n^{q,k}$ to be   nonempty is $n\ge 
q\ge 
k\ge 0.$  An immediate observation is that $B_n=\cup_{q=0}^n\cup_{k=0}^q 
B_n^{q,k}$, thus $2^n=\sum\limits_{q=0}^n\sum\limits_{k=0}^q F_n(q,k),$
where $F_n(q,k)=|B_n^{q,k}|.$ The value of $F_n(q,k)$ is defined to be zero 
whenever $n< 0.$ For more results on $F_n(q,k)$ see Monimala Nej, A. 
Satyanarayana Reddy \cite{M:S}. M.A. Nyblom 
in \cite{NY1} 
denoted $S_{r}(n)$ for the set of all binary strings 
of length $n$ without  $r$-runs of ones, where $n\in \mathbb{N}$ and $r \geq 
2$, and $T_{r}(n)=|S_{r}(n)|$. For example, if  $n=3$, $r=2$, then 
$S_{2}(3)=\{000, 101,001,100,010\}$, $T_{2}(3)=5$.
\begin{theorem}\label{thm:al1}  \hfill

\begin{enumerate}
\item \label{thm:al1:1} Let $A\in \M_n^{1,1}$, where $n\ge 4.$ Then 
$exp(A)=2(t+1),$ 
 where $t=\lfloor{\frac{m(V_1)+1}{2}}\rfloor.$ 
 And $N_n^{1,1}(2(t+1))=2 
\left[\sum\limits_{q=m(V_1)}^{n-3}F_{n-3}\left(q,m(V_1)\right) 
\right].$ 

 \item \label{thm:al1:2} Let $A\in \M_n^{0,1}$, where $n\ge 4.$ Suppose 
$\min(V_2)=h$  and 
 $t=\max \left\{h-1, \lfloor{\frac{m(V_1^h)+1}{2}}\rfloor\right\},$ where 
$V_1^h=V_1 \setminus \[1,h-1\].$  
 Then $exp(A)=2(t+1).$ And  $$N_n^{0,1}(2(t+1))= 2 
\left[\sum\limits_{i=0}^{t-2}\sum\limits_{k=2t-1}^{2t}\sum\limits_{q=k}^{n-i-4} 
F_{n-i-4}(q,k)+T_{2t+1}(n-t-3) \right].$$
\end{enumerate}
\end{theorem}
\d 
Proof of Part \ref{thm:al1:1}. Let $A\in \M_n^{1,1},$ then from 
Corollary~\ref{cor:expoi}\quad \quad  $exp(A:i,j) \leq 2$ for all $i,j\in V_2$ and for
 $j\in V_2$, $i\in V_1$ $exp(A:i,j)\leq (i-i')+2$  where $i'= max\{\[1,i\] \cap 
V_2\}$. Now it  is  sufficient to 
find $exp(A:i,j)$ for all $i,j\in V_1.$  
 Suppose $m(V_1) \neq 0.$ Then
 by definition of $m(V_1)$ there exists  $U\subseteq V_1$ such that 
 $|U|=m(V_1).$ Suppose $U=\{i_1,i_2,\ldots,i_{m(V_1)}\},$ where $i_a=i_{a-1}+1, 
a\in \[2,m(V_1)\].$
 The result follows by observing the fact that 
 $exp(A:i,j)\le exp(A:i_{t},i_{t})=2 t+2,$ for all $i,j\in V_1.$ 
 
 Before finding $N_n^{1,1}(2(t+1))$ one can observe that  the exponent of a matrix 
depends on  $m(V_1).$
 Since $a_{n,1}a_{n,n-1}a_{n,n}>0$, there is a bijection between sets $ 
\M_n^{1,1}$ 
and $B_{n-3}.$ 
 In particular there is a  bijection between $B_{n-3}^{q, m(V_1)}$ and
 $\{A\in  \M_n^{1,1} : exp(A)=2\lfloor{\frac{m(V_1)+1}{2}}\rfloor+2\}.$  To 
conclude one has to invoke Remark~\ref{rem:nn-1}.

 Proof of Part \ref{thm:al1:2}. Suppose $r=\lfloor \frac{m(V_1^h)+1}{2}\rfloor.$
 From the proof of first part we have $exp(A:i,j)\le 2(r+1)$  for $i,j\in 
V_1^h.$ Further there exists 
 $i\in V_1^h$ with $exp(A:i,i)=2r+2.$
 Thus the proof follows by observing  that $exp(A:j,j') \leq exp(A:1,1)=2h$  for $j, 
j' \in \[1,h-1\]$ and  $exp(A:j,j') \leq h+r+1$, where $j \in 
\[1,h-1\]$ and $j' \in V_1^h$.

Finally the number of primitive matrices with a given exponent follows easily from 
the definitions of 
$F_n(q,k)$ and $T_{r}(n).$
\dd

\ex
Consider the graphs in the Figure \ref{fig:figure4}. Both of them have exponent six.
\begin{figure}[htb]
\begin{center}
\includegraphics[scale=.8]{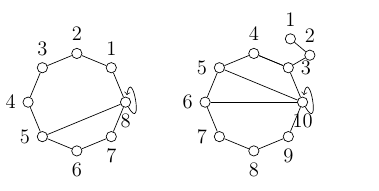}
\end{center}
\caption{Primitive symmetric companion matrices with exponent six}  \label{fig:figure4}
\end{figure}

For $n=7$, number of matrices with exponent six is given by 
$N_7^{0,1}(6)=2\left[\sum\limits_{k=3}^{4}\sum\limits_{q=k}^{3}F_{3}(q,k)+T_{5}
(2)\right]=10.$ 
Thus there are $10$ matrices in $\M_{7}^{0,1}$  with exponent $6.$
\exx

\begin{cor}
  Let $n \geq 4$. Then $E(\M_n^{1,1})=\{2(t+1) : t\in \[0, 
\lfloor\frac{n-2}{2}\rfloor\]\}$ and
  $E(\M_n^{0,1})=\{2t : t\in\[2, n-1\]\}.$
 \end{cor}
It is easy to check that the exact upper bound for the class of symmetric 
primitive matrices in \cite{shao} belongs to $E(\M_n^{0,1}).$  And the graph which achieve this  bound is in Figure \ref{fig:figure5}. 
\begin{figure}[htb]
\begin{center}
\includegraphics[]{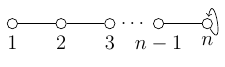}

\end{center}
\caption{The path graph of length $n-1$ with a loop at vetrex $n$}  \label{fig:figure5}
\end{figure}

%%%%%%%%%%%%%%%%%%%%%%%%
\section{Exponents of matrices in $\M_n^{1,0}$}\label{sec:CBNL}
%%%%%%%%%%%%%%%%%%%%%%%%%% 

If  $A \in \M_n^{1,0},$ then we know that $a_{n,n-1}$ is either $1$ or 
$2.$
 Thus the last row of $A$ is of the form  $1Ya_{n,n-1}0$, where $Y\in 
B_{n-3}.$ Recall that we denote an element in $\M_n^{1,0}$ as $A_Y.$  
From  Corollary \ref{cor:even} we know that if $A\in\M_n^{1,0}$, then 
$exp(A)=2d$ for some $d\in \N.$
In particular, we have the following result.
\begin{theorem}\label{thm:10expset}
 Let $n\in \N, n\ge 4.$ Then $E(\M_n^{1,0})=\{2t : t\in \[1, 
\lfloor\frac{n-1}{2}\rfloor\]\}.$
\end{theorem}
\d If $Y\in B_{n-3}^{k,k}$ and $k$ is even, then we claim that 
$exp(A_Y)$ is $k+2$ or equivalently 
 $m(V_1)+2.$
 For example, if $Y=\underbrace{11\cdots 11}_{n-3} \in B_{n-3}^{0,0}$, then  
 $A_Y=\begin{bmatrix}
        0 &1 &0 &0& \dots & 0 &1\\
        1&0 &1 &0& \dots & 0 &1\\
        0 &1 &0 &1& \dots & 0 &1\\
        \vdots &\vdots &\vdots &\vdots &\dots & \vdots &\vdots\\
        0 &0 &0 &0 &\dots & 0 &2\\
        1 &1 &1 &1&\dots & 2 &0     
    \end{bmatrix}$ and clearly  $A_Y^2>0.$ Hence $2\in E(\M_n^{1,0}).$ We now 
suppose that
$k\in \[2,n-3\]^e$.
    In this case $V_1=\[a+1,a+k\]$ for some $1\le a\le n-k-2$ and 
    $V_2=\[1,a\]\cup \[a+k+1,n-1\]$. As a consequence  if 
$i,j\in V_2$ then $exp(A_Y:i,j)\le 2.$
   Furthermore, it is easy to see that  for every $i,j\in \[1,n\],$  $exp(A_Y:i,j)\le 
exp(A_Y: p,p)=k+2$ where $p=a+\frac{k}{2}+1$ or $a+\frac{k}{2}.$ Finally the 
upper bound is achieved from   Corollary $1$ in Delorme and Sole 
\cite{D:S}. Hence the result follows.
   \dd

Next result is the   main result of this section, which 
expresses  the exponent 
of $A\in \M_n^{1,0}$ in terms of 
certain parameters related to $V_1$. It can be used as an algorithm to 
determine 
the exponent of a primitive matrix of this class.
To define these parameters, recall that we expressed $V_1$ as 
$V_1=U_1\cup U_2\cup \dots \cup U_r,$ 
where $U_i$ is a set of  consecutive numbers in $V_1$ and 
for  $i > j\in \[1, r\],$  $ |min(U_i)-max(U_j)|\ge 2.$ 
Then the maximum odd number among $|U_1|,|U_2|, \ldots, |U_r|$ is denoted  
$mo(V_1).$ Thus $mo(V_1)=m(V_1)$, if $m(V_1)$ is odd. 
Because of the  primitivity of $A,$  $ |U_i|$ is even  for some $i.$

Also,  let $q_{V_1}= $ the smallest even number in 
$\{|U_1|,|U_2|, \ldots, |U_r|\}$, and let $se(V_1)=\begin{cases}
                                  0 \;\;\; \hspace*{5mm}\mbox{if $i,i+1\in V_2$ 
for some $i\in 
\[1,n-2\],$} \\ q_{V_1}\hspace*{5mm} \mbox{otherwise.}                          
 \end{cases}$ \\
That is  $se(V_1)+3$ is the smallest odd cycle length  in $D(A).$  Now 
corresponding to each $U_i\;\;1\le i\le r$
 there is an elementary cycle  denoted  $c_{U_i}$ (or simply $c_i$),  with its 
vertex set  $V(c_{U_i})=U_i  \cup \{n,p_i,q_i\}$, 
 where $p_i,q_i\in V_2.$ If $p_i<q_i$, then we call $p_i$ and $q_i$ as the first 
vertex and the last vertex of $c_{U_i}$ respectively.

It is clear that $\ell(c_{U_i})$, the length of the cycle $c_{U_i}$  is equal to 
$|U_i|+3.$ Let $C(A)=\{c_{U_t} : 1\le t \le r\} \cup \{\text{cycles of length 
$3$}\}.$ Since there exists 
$i\in \[1,r\]$ such that $|U_i|$ is even, hence
there exists $c\in C(A)$ such that $\ell(c)$ is odd.  From this observation it 
is easy to see, for  each cycle 
$d_1\in C(A)$ with an even length  there exists 
a cycle $d_2\in C(A)$ with $\ell(d_2)=se(V_1)+3$ such that  $\ell(d_1+d_2)$ is 
$\ell(d_1)+se(V_1)+1$ or $\ell(d_1)+se(V_1)+3,$ where $d_1+d_2$ is a cycle in 
$D(A)$ with
$V(d_1+d_2)=V(d_1)\cup V(d_2).$ Now  $V(d_1+d_2)$ contains at least three 
vertices from $V_2.$ In particular, $|V(d_1+d_2)\cap V_2|=3$ or $4.$ Hence we 
let the first vertex and last 
vertex of $d_1+d_2$ be  same as the  first vertex and last vertex   of $d_1$ 
respectively. A cycle  $d_2$ with this property is called   an {\em 
associate} of 
$d_1.$ Let $C'(A)=\{c\in C(A) : \ell(c) \mbox{ is odd}\}\cup 
\{d_1+d_2 : d_1,d_2\in 
C(A), \ell(d_1) \;\mbox{ is even}, d_2 \mbox{ is an associate of $d_1$}\}.$

\ex
Consider the graph in Figure \ref{fig:figure6}. 
\begin{figure}[htb]
\begin{center}
\includegraphics[scale=.8]{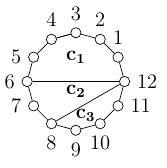}
\end{center}
\caption{An example illustrate the parameters}  \label{fig:figure6}
\end{figure}

Here $V_1=\{2,3,4,5\} \cup \{7\}\cup\{9,10\}$ and $C(A)=\{c_1, c_2, c_3\}$ with 
$V(c_1)=\{1,2,3,4,5,6,12\}$, 
$V(c_2)=\{6, 7, 8, 12\}$, $V(c_3)=\{8, 9, 10, 11,12\}$. Hence $\ell(c_1)=7$ and 
$1$, $6$ are the first vertex,  last vertex of  $c_1$  respectively.  Also 
$m(V_1)=4$, 
$mo(V_1)=1$, 
$se(V_1)=2.$ Hence $C'(A)=\{c_1,c_2+c_3,c_3\}.$
\exx

%%%%%%%%%%%%%%%%%%%%%%%%%%
\begin{theorem}\label{thm:Main:Mn10}
Let $n \geq 4$ be an integer and $A\in \M_n^{1,0}.$
\begin{enumerate}
 \item If $m(V_1)$ is odd,  then  
 $$exp(A)=\begin{cases}
 m(V_1)+se(V_1)+5 & \begin{array}{l}\mbox{if there exists $ c \in C'(A)$ with}\\ \,\,\,\, c=d_1+d_2, \ell(d_1)=m(V_1)+3\end{array}\\
 & \quad \mbox{and $|V(d_1+d_2)\cap V_2|=4$,}\\ 
 m(V_1)+se(V_1)+3 & \mbox{otherwise.}                              
 \end{cases}$$
 \item If $m(V_1)$ is even, and $mo(V_1)\ge m(V_1)-se(V_1)-1$, then 
$$exp(A)=\begin{cases}
 mo(V_1)+se(V_1)+5 &
 \begin{array}{l}\mbox{if there exists $ c \in C'(A)$ with}\\ \,\,\,\, c=d_1+d_2, \ell(d_1)=mo(V_1)+3\end{array}\\
 & \quad \mbox{and $|V(d_1+d_2)\cap V_2|=4$,}\\
 mo(V_1)+se(V_1)+3 & \mbox{otherwise.}                              
 \end{cases}$$ 
 \item In the remaining cases $exp(A)=m(V_1)+2.$
 \end{enumerate}
\end{theorem}
\d Recall that  if  $A$ is a symmetric primitive  matrix     and  
the $ij$-th entry of $A^s$ is positive , then  the $ij$-th entry of $A^{s+2}$ is positive.  Further if $i,j\in V(c)$, where $c$ 
is a cycle of 
odd  length, then there are two paths from $i$ to $j$ in $c$ such that   sum of their lengths is  $\ell(c)$.  
Hence from Corollary~\ref{cor:expoi} the smallest 
odd cycle containing 
$i,j$ plays an important 
role in finding $exp(A:i,j).$

We suppose that $c \in C'(A)$ with  $p$ as its  first vertex. Then for all $i, j 
\in V(c)$, we have

$$ exp(A:i,j) \leq \begin{cases}
                                  \ell(c)-1 & \mbox{when}\; c\in C(A),\\
                                  \ell(c)+se(V_1) & \mbox{when}\begin{array}{l} c=d_1+d_2, \mbox{ and } \\
\ell(c)=\ell(d_1)+se(V_1)+1,\end{array}\\      
     \ell(c)+se(V_1) +2& \mbox{ otherwise.}                                   
             \end{cases}$$ 
And from 
Proposition 
~\ref{pro:expoi}, equality holds  when $i=j=p+\lfloor 
\frac{\ell(c)-1}{2}\rfloor.$ Now we 
attempt to find a similar bound for $exp(A:i,j)$ when 
$i$,$j$ belongs to different cycles in $C'(A).$

 To prove the result it is sufficient to show that if  $\ga, \de\in 
C'(A)$ with $\ell(\ga) \geq \ell (\de)$, then $exp(A:i,j) \leq \ell(\ga)-1$, 
for 
$i \in V(\ga)$ and $j \in V(\de)$.

Let $p_\ga$ and $p_\de$ be the first vertex of $\ga$ and $\de$ respectively. 
Suppose $q,q'$ are nonnegative integers such that $\ell(\ga)-q-2>q$ and 
$\ell(\de)-q'-2>q'$.
For $q \leq q'$, from vertex $p_\ga+q$ to $p_\de+q'$ there are paths of length 
$q+q'+2$ and $\ell(\ga)+q'-q$ which are in opposite parity because 
$\ell(\ga)-2$ 
is 
odd. Hence $exp(A:p_\ga+q,p_\de+q') \leq max\{q+q'+2, \ell(\ga)+q'-q\}-1 \leq 
\ell(\ga)-1.$
For $q > q'$ there are paths of length $q+q'+2$ and $\ell(\de)+q-q'$ from 
vertex 
$p_\ga+q$ to $p_\de+q'$ and $exp(A:p_\ga+q,p_\de+x') \leq max\{q+q'+2, 
\ell(\de)+q-q'\}-1 \leq \ell(\ga)-1.$

Thus $exp(A:i,j) \leq \ell(\ga)-1$, for $i \in V(\ga)$ and $j \in V(\de)$.  
Hence the result follows from Lemma~\ref{lemma:exp(A)}.
\dd       

It is easy to see that Theorem~\ref{thm:10expset} also follows as a corollary 
of the above theorem. In the remaining part  of  this section
\begin{enumerate}
 \item we will  evaluate $N_n^{1,0}(b)$  when $b=2$ and $b=max\[2,n-1\]^e$.
 \item And for the remaining values of $b$, we will provide bounds for 
$N_n^{1,0}(b).$
\end{enumerate}
 If $Y \in B_{n-3}^{k,k}$, then $exp(A_Y)$ is either $k+2$ or $k+3$ depending 
on 
the parity of $k.$ From \cite{M:S} it is 
known that for $k \geq 1$
 $F_{n}^{k,k}=(n-k)+1.$ 
 Hence we have 
$4(n-b)+2\le N_n^{1,0}(b).$ And $N_n^{1,0}(b) \leq 2 \left[
\sum\limits_{k=\frac{b-4}{2}}^{b-2}\sum\limits_{q=k}^{n-3}F_{n-3}^{q,k}\right]$ 
follows 
from Theorem 
~\ref{thm:Main:Mn10}. In order to improve the upper  bound of $N_n^{1,0}(b)$  
we 
partition the set 
$\M_n^{1,0}$ as $\M_n^{1,0}=\cup_{k=0}^{n-3} \M_{n,k}^{1,0},$ where 
$\M_{n,k}^{1,0}=\{A_Y 
\in \M_n^{1,0} : Y \in B_{n-3}^{q,k}\}.$

Note that if  $A \in \M_{n,k}^{1,0}$, 
then $m(V_1)=k$ and $se(V_1)\in \[0,k\]^e.$

 If $n\ge 4$, then it is easy to verify that  $2\notin E(\M_{n,k}^{1,0})$ if and 
only if 
$k\ge 1$ and  $E(\M_{n,0}^{1,0})=\{2\}.$
 As a consequence we get $|\M_{n,0}^{1,0}|=2$ and  $N_n^{1,0}(2)=2.$ 
The following result evaluates $E(\M_{n,k}^{1,0})$ for various values of $n$ 
and $k$ through which the upper bound of $N_n^{1,0}(b), b \geq 4$ can be improved.
\begin{theorem}\label{counting1}
Let $n \geq 5$ and $k$ be positive integers.
\begin{enumerate}
\item\label{general}
If $n$ is odd and $n \geq 2(k+2)+1$, or if $n$ is even and $n  \geq 3(k+2)-1$, then 
$E(\M_{n,k}^{1,0})= \[k+2,2(k+2)\]^e$.
\item\label{nk odd}
If $k$ is odd and $n \in \[k+3,2(k+2)-1\]^o$, then $E(\M_{n,k}^{1,0})= 
\[k+2,n-1\]^e$.  
\item\label{k ev n od}
If $k$ is even and $n \in \[k+3,2(k+2)-1\]^o \setminus \{k+5, k+7\}$, then $ 
E(\M_{n,k}^{1,0}) = \[k+2,n-1\]^e\setminus \{k+4\}$. And for $n \in 
\{k+5, k+7\}$, $ 
E(\M_{n,k}^{1,0}) = \[k+2,n-1\]^e.$
\item\label{k od n ev}
If $k$ is odd and $n \in \[k+5,3(k+1)\]^e$, then $E(\M_{n,k}^{1,0}) = 
\[k+3,l\]^e$,  \\
where $l=\begin{cases}
 k+2\lfloor \frac{n-k-4}{4} \rfloor+3 & \mbox{if} \;n-k-4=4\lfloor 
\frac{n-k-4}{4} \rfloor+1,\\
 k+2\lfloor \frac{n-k-4}{4} \rfloor+5 & \mbox{if}\; n-k-4=4\lfloor 
\frac{n-k-4}{4} \rfloor+3. 
                                                               \end{cases}$

\item\label{nk ev}
If $k$ is even and  $n \in \[k+4, 3k+2\]^e$, then $E(\M_{n,k}^{1,0}) = 
\[k+2, n-k\]^e$  \; but if $n-k \leq k+2$, then $E(\M_{n,k}^{1,0}) = \{k+2\}.$ \\
For $k$ is even and $n=3(k+1)+1$, $E(\M_{n,k}^{1,0}) = \[k+2, 2(k+1)\]^e.$
\end{enumerate}
\end{theorem}

\d
First  observe that the  lower bound of 
$E(\M_{n,k}^{1,0})$ will be 
attained by $A_Y$  where $Y\in B_{n-3}^{k,k}.$

Proof of part ~\ref{general}. For any element in $\M_{n,k}^{1,0}$ the 
maximum value of 
 $se(V_1)$  is either  $k$ or $k-1$ depending on whether  $k$ 
is even or odd respectively. Then from Theorem ~\ref{thm:Main:Mn10} the upper bound  
$2(k+2)$ belongs to $E(\M_{n,k}^{1,0})$ whenever $n \geq 2(k+2)+1$ when $n$ is odd 
or $n \geq 3(k+2)-1$ when $n$ is even.
\begin{itemize}
 \item Suppose $n \geq 2(k+2)+1$ is odd and  $$Y=(\underbrace{0, \dots,0}_{k},1, 
0, 1,\underbrace{0,\dots,0}_{t},\underbrace{1, 0,\dots,1,0}_{n-k-t-6}),$$ where $t 
\in \[0,k-1\]$ and  $2\nmid t+k$. Then the exponent of  $A_Y \in 
\M_{n,k}^{1,0}$  is   $k+t+5.$
\item Suppose $n \geq 3(k+2)-1$ is even and $k$ is odd. If we choose 
$$Y=(\underbrace{0, 
\dots,0}_{k},1,0,1,\underbrace{0,\dots,0}_{t},1,\underbrace{0, 
\dots,0}_{k-1},1,\underbrace{0,1,\dots, 0}_{n-2k-t-7}),$$ where $t \in \[0, 
k-1\]^e$, then the exponent  $A_Y$ is  $k+t+5$. Similarly for $k$ is 
even,  the exponent $A_Y$ is  $k+t+5$, whenever
$$Y=(\underbrace{0, \dots,0}_{k},1,\underbrace{0, 
\dots,0}_{k},1,0,1,\underbrace{0,\dots,0}_{t},\underbrace{1, 
0,\dots,1,0}_{n-2k-t-7})$$ and $t \in \[1, k-1\]^o.$ 
\end{itemize}
Finally,  $A_Y\in \M_{n,k}^{1,0}$ has exponent $k+4$, if 
$$Y=(\underbrace{0, \dots,0}_{k-1},1, \underbrace{0,\dots,0}_{k},\underbrace{1, 
1,\dots,1}_{n-2k-3}).$$

Proof of part ~\ref{nk odd}. Suppose $n \neq k+4$ and $$Y=(\underbrace{0, 
\dots,0}_{k},1, 
0, 1,\underbrace{0,\dots,0}_{t},\underbrace{1, 0,\dots,1,0}_{n-k-t-6}),$$ where $t 
\in \[0,n-k-6\]^e.$  Then $exp(A_Y)$ is $k+t+5$ and $E(\M_{n,k}^{1,0})= 
\[k+2,n-1\]^e.$

Proof of part ~\ref{k ev n od}. For $n \in \[k+9, 
2(k+2)-1\]^o$, there are no elements in $\M_{n,k}^{1,0}$ with exponent 
$k+4.$ 
Because then  $mo(V_1)+se(V_1)=k+1$ or $mo(V_1)+se(V_1)=k-1$. And 
$mo(V_1)+se(V_1)=k+1$ implies $se(V_1)=k$ and $mo(V_1)=1$ and the exponent of the 
corresponding elements is $k+6$. Since $n-k-4 \leq k-1$, 
$mo(V_1)+se(V_1)=k-1$ is not possible. For $n=k+5$ and $n=k+7$, $A_Y$ has 
exponent $k+4$, where $Y=(\underbrace{0, \dots,0}_{k},1, 0) \in B_{k+2}$ and 
$Y=(0,1,\underbrace{0, \dots,0}_{k},1,0) \in B_{k+4}$ respectively. And for $n 
\in \[k+7, 2(k+2)-1\]^o$, $exp(A_Y)=k+t+5$, where $Y=(\underbrace{0, 
\dots,0}_{k},1, 
0, 
1,\underbrace{0,\dots,0}_{t},\underbrace{1, 0,\dots,1,0}_{n-k-t-6})$ and $t 
\in \[0,n-k-6\]^o.$ Hence the result follows in this case.

Proof of part ~\ref{k od n ev}. Suppose $d=n-k-4$ and 
$d'=\lfloor\frac{d}{4}\rfloor$. Then for this part either $d=4d'+1$ or $d=4d'+3$ 
and $se(V_1) \in \[0, 2d'\]^e$.

For $d=4d'+1,$ consider $Y=(\underbrace{0, 
\dots,0}_{k},1,\underbrace{0,\dots,0}_{t},1,\underbrace{0, 
\dots,0}_{2d'},\underbrace{1,0,1,\dots,0}_{n-k-t-2d'-5})$, where $t \in \[0, 
2d'\]^e$. Then $exp(A_Y) = k+t+3$. If $d=4d'+3$, then the exponent of $A_Y$ is 
$k+t+5$, where 
$$Y=(\underbrace{0, 
\dots,0}_{k},1,0,1,\underbrace{0,\dots,0}_{t},1,\underbrace{0, 
\dots,0}_{2d'},\underbrace{1,0,1,\dots, 0}_{n-k-t-2d'-7})$$ and $t \in \[0, 
2d'\]^e$. Hence the result follows in this case.

Proof of part ~\ref{nk ev}. Let $n \in \[k+4,2(k+1)\]^e$. Then for the elements 
in 
$\M_{n,k}^{1,0}$ either $mo(V_1)$ does not exists or 
$mo(V_1)+se(V_1)\leq 
m(V_1)-3$, which gives $E(\M_{n,k}^{1,0})=\{k+2\}$.

For $n \in \[2(k+2), 3k+4\]^e$, maximum value of $mo(V_1)+se(V_1)$ is $n-k-5$. 
Then the upper bound for the set $E(\M_{n,k}^{1,0})$ is $n-k$ or $n-k-2$ 
according as $n \in \[2(k+2), 
3k+2\]^e$ or $n=3k+4$ respectively. And this upper bound will be attained by $A_Y$, 
where 
$Y=(\underbrace{0, \dots,0}_{k-1},1, 
\underbrace{0,\dots,0}_{k},1,\underbrace{0,\dots,0}_{n-2k-4})$. Also for 
$2(k+2) 
\leq n \leq 3k+4$, exponent of $A_Y$ is $k+t+5$, where 
$Y=(\underbrace{0, \dots,0}_{k},1,\underbrace{0, 
\dots,0}_{k},1,0,1,\underbrace{0,\dots,0}_{t},\underbrace{1, 
0,\dots,1,0}_{n-2k-t-7})$ with $t \in \[1, n-2k-7\]^o$. And the exponent of $A_Y$ is 
$k+4$, where $Y=(\underbrace{0, 
\dots,0}_{k-1},1, \underbrace{0,\dots,0}_{k},\underbrace{1, 
1,\dots,1}_{n-2k-3}).$ Hence the proof is complete.
\dd
Recall that $N_n^{1,0}(b)\le 2 \left[
\sum\limits_{k=\frac{b-4}{2}}^{b-2}\sum\limits_{q=k}^{n-3}F_{n-3}^{q,k}\right].$ 
Using the above theorem,  by restricting some of the 
values of $k$ in the first summation, we will improve this upper bound of $N_n^{1,0}(b)$ for 
certain values of $b$.  

Let $S_n^{1,0}(b)$ denote the set of all 
possible values of $k$ in the first  summation.  Then $S_n^{1,0}(b)\subseteq 
\[\frac{b-4}{2}, b-2\].$ The following result provides all possible values in 
$S_n^{1,0}(b).$

\begin{cor} \label{cor:UB}
Let $b \in E(\M_n^{1,0})\setminus \{2\}$. 
\begin{enumerate}
\item  \label{cor:UB:1} Then $b-2\in S_n^{1,0}(b).$
 \item \label{cor:UB:2} If $n \in \[b+5, 2b-5\]^o$, then $S_n^{1,0}(b)= 
\[\frac{b-4}{2}, b-2\] 
\setminus \{b-4\}$.  
\item \label{cor:UB:3}
Let $n$ be even. Then by Part ~\ref{general} of the above theorem $\[\frac{b-4}{2},  
\lfloor \frac{n+1}{3}\rfloor-2\]\subseteq S_n^{1,0}(b).$
For the remaining cases 
that 
is when $ k \geq \lfloor \frac{n+1}{3}\rfloor-1$, we have the following.
\begin{enumerate}
\item
If $k$ is odd,  then $k\in S_n^{1,0}(b)$,  provided  $l \geq b$, where $l$ is  as defined  in Part ~\ref{k od n ev} of the above theorem.
\item 
If $k$ is even, $k \leq n-b$,  and $n \in \[2k+4, 3k+2\]^e$,  then $k\in S_n^{1,0}(b)$.    If $n=3k+4$ and $b \leq \frac{2(n-1)}{3}$, then $k \in S_n^{1,0}(b).$
\end{enumerate}
\end{enumerate}
\end{cor}

For example,   $S_{20}^{1,0}(12)=\[4,10\].$ Hence there is  no improvement,  where 
as $S_{20}^{1,0}(16)=\{9,11,13,14\},$ a significant improvement.

%%%%%%%%%%%%%%%%%%%%%%%%%%%%%%%%%%%%%%%%%%%%%%%%
\begin{theorem}
Let $n \in \N$ and $b=max\[2,n-1\]^e$. Then $$N_n^{1,0}(b)=\begin{cases}
                
                   2(2n-7) & \mbox {if $n \geq 5$, \text{$n$ is odd}}, \\
                 18 & \mbox {if $n \geq 8$, \text{$n$ is even.}}

                \end{cases}$$
\end{theorem} 
\d  We divide the proof into the even and odd cses.
\begin{description}
 \item[ $n$ is odd.]
 From Part~\ref{cor:UB:2} of  Corollary~\ref{cor:UB},  a necessary condition 
for 
a matrix $A_Y\in 
\M_{n,k}^{1,0}$ having  exponent 
 $n-1$  is $k\in  \[\frac{n-5}{2}, n-3\]$ where $Y\in 
B_{n-3}^{q,k}.$ 
 Here $q \leq 2(k+1)$ follows from the fact that  $k \leq q \leq n-3$ and  
$n-k-4 \leq k+1.$

Since exponent of $A_Y$ is $n-1$, hence from  Theorem \ref{thm:Main:Mn10}  we 
have
\begin{itemize}
 \item if $k\in \[\frac{n-3}{2}, n-6\]$, then $q\in 
\{n-4,n-5\}.$
\item if  
$k=\frac{n-5}{2}$, then  $q=n-5$. 
\item if  $k \in\{n-5, n-4\}$, then  $q=n-4$ and when  $k=n-3$, $q$ is also  
$n-3$.
\end{itemize}

 Thus $A_Y$ will have exponent $n-1$ only if $Y\in S,$  where $$S=\bigcup 
\limits_{k=\frac{n-5}{2}}^{n-6}B_{n-3}^{n-5,k} 
\bigcup \limits_{k=\frac{n-3}{2}}^{n-4}B_{n-3}^{n-4,k} \bigcup 
B_{n-3}^{n-3,n-3}.$$

When $q=k=n-3$ there are only two  matrices with exponent 
$n-1,$ in the remaining cases  there are four matrices with exponent $n-1$.
Hence the result follows in this case.

\item[ $n$ is even.]
Suppose  $n \geq 10$.  From Part~\ref{cor:UB:3} of  Corollary~\ref{cor:UB} 
a necessary condition 
for a matrix $A\in \M_{n,k}^{1,0}$ having  exponent 
$n-2$  is  $k \in \{n-7,\; n-5, \;n-4\}$. 
\begin{itemize}
 \item When  $k= n-4, \; n-5$ it is easy to see that $q=n-4$ and $q=n-5$ 
respectively 
and exponent of $A_Y$ is $n-2$ for all  $Y \in B_{n-3}^{k,k}$. 
\item When  $k=n-7$, we have $n-7 \leq q \leq n-4$. But $exp(A_Y)\ne n-2$ if 
$q=n-7$ or $q=n-4.$
Further  when  $q=n-6$  or $q= n-5$, there are elements in 
$\M_{n,n-7}^{1,0}$ whose exponent is $n-2$ 
provided  $se(V_1)=0.$
\end{itemize}

 Thus $A_Y \in 
\M_{n}^{1,0}$  will have exponent $n-2$ only if $Y\in S,$ where 
$$S=\bigcup \limits_{q=n-6}^{n-5}B_{n-3}^{q,n-7} \bigcup B_{n-3}^{n-5,n-5} 
\bigcup B_{n-3}^{n-4,n-4}.$$
For $n=8$ possible values of $q$ and $k$ will be same as above except $k=1$ 
whenever $q=3$. And modified value will be $k=2$ and $q=3$. Now for $n \geq 8$ 
it is easy to check that when $k=n-7$, then for each $q$ there are four 
matrices 
with exponent $n-2.$ Hence the result follows in this case.
%%%%%%%%%%%%%%%%%%%%%%%%%%%
\end{description}
\dd 

When  $n=6$, the number of primitive  matrices with exponent $4$ is $10.$

%%%%%%%%%%%%%%%%%%%%%%%%%%%%%%%%%%%%%%%%%%%%%%%%%%%%%%%
\section{Exponents of matrices in $\M_n^{0,0}$}\label{sec:NCNL}
%%%%%%%%%%%%%%%%%%%%%%%%%
Suppose that $n \geq 4$ and $A\in \M_n^{0,0}$. Then $exp(A)\ge 4$ as  $exp(A:1,1)\ge 
4.$ Liu et al. \cite{B:M:W}  showed that the exponent of  a 
primitive symmetric matrix with trace zero is at most
$2(n-2).$ Hence  
$E(\M_n^{0,0}) \subseteq \[4, 2(n-2)\]^e.$
The following result  shows the equality of these sets. Proof of the result is 
available in \cite{B:M:W}. We 
include the proof for the sake of completeness.
\begin{theorem}
 Suppose $n\geq 4$. Then $E(\M_n^{0,0})=\{2t : t\in \[2, n-2\]\}.$
\end{theorem}
\d
 If  $A\in \M_n^{0,0},$ then $a_{n,1}=a_{n,n}=0$ and $a_{n,n-1}=1$ or 
$2.$ 
If  the  last row of $A$ is 
 $[0 \underbrace{0\, 0\dots0\, 0}_{t-1}\underbrace{1\, 1\dots1\, 1}_{n-3-(t-1)} 2 \,0],$ where 
$t\in \[1,n-3\]$, then it is easy to see that 
 $exp(A)=exp(A:1,1)=2(t+1).$
\dd

Now it is easy to see that the upper bound for the class of  symmetric primitive 
$(0, 1)$- matrices with zero trace is $2(n-2)$. And this upper bound will be 
attained by the graph in Figure \ref{fig:figure7}.
\begin{figure}[htb]
\begin{center}
\includegraphics[scale=.8]{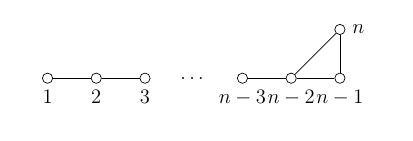}
\end{center}
\caption{A lollipop graph}  \label{fig:figure7}
\end{figure}

 Our next goal is to find $exp(A)$ for $A\in \M_n^{0,0}$, $n \geq 4$.  We 
will continue with the notation used in previous sections.
 Recall that  $\min(V_2)=h$ and $V_1^h = V_1 \setminus 
\[1,h-1\].$  It is clear that   $2 \leq h \leq n-2$. We denote the cycle in  
$C(A)$ whose first vertex is $h$ as $c^h.$
The following observations are easy to verify. 
\begin{enumerate}
\item
For all $i,j \in  \[1,h-1\]$,  $exp(A:i,j) \leq 2(h-1)+exp(A:h,h).$  The equality 
holds for $i=j=1$. Also $$exp(A:h,h)=\begin{cases}
                
                  se(V_1^h)+2 & \mbox {if $\ell(c^h)=se(V_1^h)+3$,} \\
                 se(V_1^h)+4 & \mbox {\text{otherwise.}}
                             
                \end{cases}$$

\item
For all $i, j \in V \setminus \[1,h-1\]$, $exp(A:i,j)$ can be found from the previous 
section.
\item 
For $i \in \[1,h-1\]$ and $j \in V \setminus 
\[1,h-1\]$, $exp(A:i,j) \leq (h-1)+exp(h,j)$. The equality holds for $i=1$.
\end{enumerate}

Thus to find $exp(A)$, it is sufficient to find  $exp(A:h,j)$, where $j 
\in V \setminus \[1,h\].$ It is now clear that we will 
use the results of the previous section to find $exp(A:h,j)$ but now they are 
applied to $V_1^h,$ instead of $V_1.$
\ex
For an example, consider the graph in Figure \ref{fig:figure8}. There we have $V_1=\{1,3,4,6,8,9\}$, 
$V_2=\{2,5,7,10\}$, $min(V_2)=h=2$. Hence $V_1^h=\{1,3,4,6,8,9\}\setminus 
\{1\} =\{3,4,6,8,9\}=\{3,4\} \cup \{6\} \cup \{8,9\}$ and for all $i, j \in 
V \setminus \{1,2\}=\[3, 11\]$, $exp(A:i,j) \leq 6$. Also  $V(c^2)=\{2,3,4,5,11\}$ 
be the vertex set of 
elementary 
cycle $c^2$ and $exp(A:2,2)=4$, $exp(A:1,1)=2+exp(A:2,2)=6.$ Thus to find 
$exp(A)$, it is now sufficient to find  $exp(A:2,j)$, where $j \in \[3, 11\].$ 
\exx
\begin{figure}[htb]
\begin{center}
\includegraphics[scale=.8]{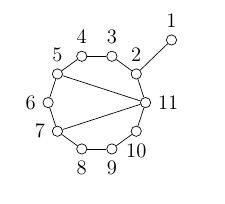}
\end{center}
\caption{An Example}  \label{fig:figure8}
\end{figure}
The following result provides an  upper bound on $exp(A:h,i)$ whenever $i$  
belongs to the vertex set of an elementary cycle. 
\pagebreak
\begin{lemma}\label{Lem:l2}
Let $A\in \M_n^{0,0}$ and $c \in C(A)$, $i \in V(c).$
\begin{enumerate}
 \item \label{lem:l21}If $\ell(c)$ is even,  then  
 $$exp(A:h,i) \leq \begin{cases}
 \frac{\ell(c)}{2}+se(V_1^h)+1 & \mbox{if $\ell(c^h)=se(V_1^h)+3,$ }\\
 \frac{\ell(c)}{2}+se(V_1^h)+3 & \mbox{\text{otherwise.} }            
                  
 \end{cases}$$
 \item \label{lem:l22}If $\ell(c)$ is odd and $\ell(c)>se(V_1^h)+3$, then 
$$exp(A:h,i) \leq \begin{cases}
  \frac{\ell(c)+se(V_1^h)+1}{2} & \mbox{if $\ell(c^h)=se(V_1^h)+3,$ }\\
 \frac{\ell(c)+se(V_1^h)+3}{2} & \mbox{\mbox{\text{otherwise.} } }            
                                              
 \end{cases}$$ 
 \item \label{lem:l23}If $\ell(c)=se(V_1^h)+3$, then 
$exp(A:h,i) \leq se(V_1^h)+2.$
 \end{enumerate}
\end{lemma}
\d 
Proof of Part~\ref{lem:l21}. Let $p$ be the first vertex of $c$ and  
$i=p+\frac{\ell(c)-2}{2}.$  
Let $h-n- p- (p+1)-(p+2)- \dots-(i-1)-i$ and  $h-n-(p+\ell(c)-2)-(p+\ell(c)-3) 
\dots-(i+1)-i$ 
be two  paths from $h$ to $i$ of   length $\frac{\ell(c)-2}{2}+2$ such that  
together  cover all the vertices of $c.$
Now our objective is to find another pair of paths from $h$ to $i$ with parity 
different from that of $\frac{\ell(c)-2}{2}+2$ and containing all the 
vertices of $c.$
If $\ell(c^h) = se(V_1^h)+3$ then there exist a path of length 
$\frac{\ell(c)-2}{2}+se(V_1^h)+3$, otherwise there is a path of 
length $\frac{\ell(c)-2}{2}+se(V_1^h)+5$ satisfying the required conditions.
Hence in this case  the result follows by applying Corollary~\ref{cor:expoi}.

Proof of Part~\ref{lem:l22}.  Suppose $\ell(c^h)=se(V_1^h)+3$ and $p$ be the 
first vertex of $c$.
Let $l_1=\frac{\ell(c)-se(V_1^h)-3}{2}$ and $l_2=\frac{l(c)+se(V_1^h)-1}{2}.$
Consider vertices  $i=p+l_1$ and $i'=p+l_2.$  
Then by a proof  similar to the proof of part~1, there are two  paths from $h$ to $i$ 
of  lengths $l_1+2$ and 
$l_1+se(V_1^h)+3$ such that first  path contains  the vertex set $\{n, 
p, p+1,p+2, 
\ldots,i\}$ and  second path contains the vertex set $\{n, 
p+\ell(c)-2, p+\ell(c)-3, \ldots,i'\}$. 

For vertex $j,$ 
$i < j < i'$, there are paths from $h$ to $j$ of lengths $r$ and $s$ such that 
$r+s=\ell(c)+2$ and $max\{r,s\}< l_1+se(V_1^h)+3 $.

For $\ell(c^h) \neq se(V_1^h)+3$, choose  $i=p+l_1-1$ and 
$i'=p+l_2+1$. And proceed in similar manner by replacing paths $l_1+2$, 
$l_1+se(V_1^h)+3$ with $l_1+1$, $l_1+se(V_1^h)+4$ respectively.
Hence from Corollary~\ref{cor:expoi} the result follows for this case.

 Proof of Part~\ref{lem:l23}. Let $i \in V(c)$.  The from 
Corollary~\ref{cor:expoi} there exists two paths from  $h$ to $i$
 of lengths $r$ and $s$ such that 
$r+s=\ell(c)+2$ 
and $max\{r,s\} \leq \ell(c).$ Hence $exp(A:h,i) \leq 
se(V_1^h)+2$ for all $i \in V(c)$.
\dd

Figure \ref{fig:figure9} illustrates the proof of part~\ref{lem:l22} of the above 
result. This graph is a subgraph of a graph with $se(V_1)=2$ and 
$\ell({c^h})=5.$ Choose $i=p+2$ and $i'=p+5$, then  there 
is a path of length $se(V_1)+1$ between  $i$ and $i'$. 
\begin{figure}[htb]
\begin{center}
\includegraphics[scale=.7]{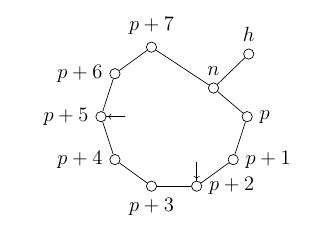}
\end{center}
\caption{A subgraph of a graph with $n$ vertices,$ \; se(V_1)=2$ and 
$\ell({c^h})=5$}  \label{fig:figure9}
\end{figure}

The following result is the main result of this section which can be used as
an algorithm to  compute $exp(A),$ where $A\in \M_n^{0,0}.$ 
\pagebreak
\begin{theorem}\label{thm:Mn001}\quad

\begin{enumerate}
 \item \label{thm:Mn001:1} Let $A\in \M_n^{0,0}$ and $m(V_1^h)$ is odd.  Then
 $$exp(A)=\begin{cases}
  2h'+se(V_1^h) & \mbox{if $h' \geq \frac{m(V_1^h)+5}{2}$,}\\
 m(V_1^h)+se(V_1^h)+5 & \mbox{if\, }\left\{\begin{array}{l} 2 \leq h' \leq \frac{m(V_1^h)+3}{2}\\\mbox{and there exists } c 
\in C'(A)\\  \mbox{with} \; c=d_1+d_2, \ell(d
_1)=\\ \quad \quad \quad \quad  \quad \quad m(V_1^h)+3
 \\ \mbox{and}\; |V(d_1+d_2) \cap V_2|=4,\end{array}\right. \\
 m(V_1^h)+se(V_1^h)+3 & \mbox{otherwise, }\\
  \end{cases}$$
 where  $h'= \begin{cases} 
              h &\mbox{ if  $\ell(c^h)=se(V_1^h)+3$,}\\
              h+1 & \mbox{ if  $\ell(c^h) \neq se(V_1^h)+3$.}
             \end{cases}$
\item \label{thm:Mn001:2} Let $A\in \M_n^{0,0}$ and $m(V_1^h)$ be even. Suppose 
$mo(V_1^h)$ does not exists 
or 
$mo(V_1^h) \leq m(V_1^h)-se(V_1^h)-3.$  Then\\
\ \\  
 $exp(A)=\begin{cases}
 m(V_1^h)+2 & \mbox{if $2 \leq h \leq \frac{m(V_1^h)-se(V_1^h)}{2},$}\\
2h+se(V_1^h) & \mbox{if } \left\{\begin{array}{l}   h \geq \frac{m(V_1^h)-se(V_1^h)+2}{2} \\ \mbox{and} \; 
\ell(c^h)=se(V_1^h)+3,\end{array}\right.
\\
 2h+se(V_1^h)+2 & \mbox{if }  \left\{\begin{array}{l}  h \geq \frac{m(V_1^h)-se(V_1^h)+2}{2}\\ \mbox{ and } 
 \ell(c^h) \neq se(V_1^h)+3.\end{array}\right.
 \end{cases}$
\end{enumerate}
\end{theorem}
\d
Proof of  Part~\ref{thm:Mn001:1}. Suppose $\ell(c^h)=se(V_1^h)+3$. Then 
$exp(A:h,h)=se(V_1^h)+2$ and
\begin{enumerate}
\item
for all $i,j \in \[1,h-1\]$, $exp(A:i,j) \leq 2h+se(V_1^h)$.
\item
for all $i,j \in V \setminus \[1,h-1\]$
$exp(A:i,j) \leq$ $$ \begin{cases}
                              m(V_1^h)+se(V_1^h)+5 & \mbox{if }  \left\{\begin{array}{l} \mbox{there exists }
c \in C'(A)\\ \mbox{with } c=d_1+d_2, \ell(d_1)=m(V_1^h)+3\\ \mbox{ and }|V(d_1+d_2) \cap V_2|=4\end{array}\right.
                 \\
                 m(V_1^h)+se(V_1^h)+3 & \mbox {otherwise.}
                 \end{cases}$$
\item
for $i \in \[1,h-1\]$ and $j 
\in V \setminus \[1,h-1\]$, $exp(A:i,j) \leq h+\frac{m(V_1^h)+3}{2}+se(V_1^h).$ 
\end{enumerate}

Now we have  $2h+se(V_1^h) \leq 
 m(V_1^h)+se(V_1^h)+3$ and $h+\frac{m(V_1^h)+3}{2}+se(V_1^h) \leq 
 m(V_1^h)+se(V_1^h)+3$ whenever $h \leq \frac{m(V_1^h)+3}{2}$. For $ h \geq 
\frac{m(V_1^h)+5}{2}$, $m(V_1^h)+se(V_1^h)+5\leq 
2h+se(V_1^h)$ and $h+\frac{m(V_1^h)+3}{2}+se(V_1^h) \leq 2h+se(V_1^h)$. Hence the 
result follows from  Lemma ~\ref{lemma:exp(A)}. For $\ell(c^h) \neq 
se(V_1^h)+3$, 
$exp(A:h,h)=se(V_1^h)+4$ and we imitate the proof of the case  $\ell(c^h) = 
se(V_1^h)+3$ with $exp(A:h,h)=se(V_1^h)+4.$

Proof of Part~\ref{thm:Mn001:2}. Suppose $\ell(c^h)=se(V_1^h)+3$. Then 
$exp(A:h,h)=se(V_1^h)+2$. If $mo(V_1^h)$ 
does not exist or $mo(V_1^h) \leq m(V_1^h)-se(V_1^h)-3$, then 
\begin{enumerate}
\item
for all $i,j \in \[1,h-1\]$, $exp(A:i,j) \leq 2h+se(V_1^h)$.
\item
for all $i,j \in V \setminus \[1,h-1\]$, $exp(A:i,j) \leq m(V_1^h) +2$.
\item
for $i \in \[1,h-1\]$ and $j 
\in V \setminus \[1,h-1\], exp(A:i,j) \leq  h+\frac{m(V_1^h)+se(V_1^h)}{2}+1$.
\end{enumerate}

Now we have $h+\frac{m(V_1^h)+se(V_1^h)}{2}+1 \leq m(V_1^h)+2$ and 
$2h+se(V_1^h) 
\leq m(V_1^h)+2$ whenever $h \leq \frac{m(V_1^h)-se(V_1^h)}{2}$. For $ h \geq 
\frac{m(V_1^h)-se(V_1^h)+2}{2}$,   $h+\frac{m(V_1^h)+se(V_1^h)}{2}+1 
\leq 2h+se(V_1^h)$ and $m(V_1^h)+2 \leq 2h+se(V_1^h)$. Hence the result follows 
from  Lemma ~\ref{lemma:exp(A)}. For $\ell(c^h) \neq se(V_1^h)+3$, 
$exp(A:h,h)=se(V_1^h)+4$ and the rest of the proof 
is similar to the previous case.
\dd

\begin{note}\label{note1}
If $mo(V_1^h)$ exist and $mo(V_1^h) > m(V_1^h)-se(V_1^h)-3 $, then $exp(A)$ 
will 
be determined by Part~\ref{thm:Mn001:1} of Theorem ~\ref{thm:Mn001} with 
$m(V_1^h)=mo(V_1^h).$
\end{note}

The following theorem is the last result of this section. This will 
evaluate $N_n^{0,0}(4)$  and $N_n^{0,0}(2n-4 )$. Part \ref{count: upper} of the Theorem \ref{counting3} can be verified with Theorem 2 in \cite{B:M:W}.
\begin{theorem}\label{counting3}\quad
 \begin{enumerate}
  \item \label{count: upper} If  $n\ge 4,$ then $N_n^{0,0}(2n-4)=2.$
  \item \label{count: lower} If   $n\ge 5$, then \;  $2 
\left[\sum\limits_{q=0}^{\lfloor \frac{n+1}{3}\rfloor} \binom{n-2q-4}{q} 
\right]\leq N_n^{0,0}(4) \leq  \\
2 \left[\sum\limits_{q=1}^{\lfloor \frac{2n-7}{5} 
\rfloor} F_{n-5}(q,1)+\sum\limits_{q=2}^{\lfloor \frac{2n+2}{3}\rfloor} 
F_{n-5}(q,2)+1 \right].$
   
 \end{enumerate}
\end{theorem}
\d
Proof of Part ~\ref{count: upper}. From Theorem ~\ref{thm:Mn001} it is clear 
that only exponent involving  $h$ may be $2(n-2)$, because the maximum value that $m(V_1^h)+se(V_1^h)$ can be is $n-3.$ Also 
we have $2 \leq h \leq n-2$. Suppose $h=n-k$ for some $k \geq 
2$. For $k >2$ it can easily be checked that no suitable values for $se(V_1)$ exist 
which will provide matrices with exponent $2(n-2).$ For $k=2$, \; $se(V_1)=0$ 
and the corresponding matrices have exponent $2(n-2).$

Proof of Part ~\ref{count: lower}. Suppose $A\in \M_n^{0,0}$ and $exp(A)=4.$ 
Then $A$ is in one of the following classes 
\begin{enumerate}
\item
from Part~\ref{thm:Mn001:1} of Theorem ~\ref{thm:Mn001},  $h=2$, $m(V_1^h)=1$, 
$se(V_1^h)=0$, 
$\ell(c^h)=3$ 
and $|V(d_1+d_2) \cap V_2|=3$ for each cycle $c \in C'(A)$ with $c=d_1+d_2$, 
$\ell(d_1)=4$. That 
is the 
last row of $A$ is 
of form $011Y10$, where $Y \in 
\bigcup\limits_{q=1}^{\lfloor 
\frac{2n-7}{5} 
\rfloor} B_{n-5}^{q,1}$ such that each zero in $Y$ has two 
consecutive 
ones to its left or right. Here the  upper bound follows because for each $q$ there 
will be 
at least $\lfloor \frac{3x-2}{2}\rfloor$ ones.
\item
from Part~\ref{thm:Mn001:2} of Theorem ~\ref{thm:Mn001}, $h=2$, $m(V_1^h)=0$ or 
$m(V_1^h)=2$, $mo(V_1^h)$ 
does not exists, 
$se(V_1^h)=0$, $\ell(c^h)=3$. That is the  last row of $A$ is of form 
$011Y10$, where $Y \in B_{n-5}^{0,0}$ or $Y \in 
\bigcup\limits_{q=1}^{\lfloor 
\frac{n+1}{3}\rfloor} B_{n-5}^{2q,2}$. More precisely there are 
$\sum\limits_{q=0}^{\lfloor \frac{n+1}{3}\rfloor} \binom{n-2q-4}{q}$ matrices. 
\item
from Note ~\ref{note1}, $h=2$, $m(V_1^h)=2$, $mo(V_1^h)=1$, $se(V_1^h)=0$,
$\ell(c^h)=3$ and $|V(d_1+d_2) \cap 
V_2|=3$ for each cycle $c \in C'(A)$ with $c=d_1+d_2$, $\ell(d_1)=4$. That is the 
last 
row of $A$ is of form $011Y10$, where $Y \in 
\bigcup\limits_{q=3}^{\lfloor 
\frac{2n+2}{3}\rfloor} B_{n-5}^{q,2}$ and $Y$ must contains a 
subword 
of zeros 
of length $1$ and each such subword has two consecutive ones to its left or 
right.
\end{enumerate}
Hence the result follows.
\dd

%\section*{Conclusion}
{\bf Conclusion:} For a given $n,$ it is quite difficult to find out the exact value of $N_n^{\al, 0}(b), \; \al \in \{0,1\}$, $b \geq 4$ and $b \neq \max \bigg(E(\M_n^{\al,0})\bigg).$
%\newpage
%\section{Number of primitive matrices with given exponent}

\begin {table}[]
\caption {Number of Matrices $N_n^{\al, \be}(b)$ for small $n$} \label{tab:table:1}

\begin{center}

\begin{tabular}{|l|l|l|l|l|l|l|l|l|l|l|l|l|l|}
%\begin{tabular}{|l|l|l|l|l|l|l|l|l|l|l|l|l|l|l|}
\hline
&&\multicolumn{9}{|c|}{Exponent($b$)} &{\hskip-.1mm$\begin{array}{cc}\hbox{Range of $N_{n}^{0,0}(4)$}\\ \hbox{from  Thm ~\ref{counting3}} \end{array}$}
\\
\hline
$N_n^{\al, \be}(b)$& $\al,\be$&2\hspace*{2mm}&4\hspace*{2mm}&6\hspace*{2mm}&8\hspace*{2mm}
&10&12&14&16&18&\\
\hline
&\multirow{2}{*}{}$1,1$ &2&&&&&&&&& \\
 & $0,1$ &&2&&&&&&&& \\
$N_3^{\al, \be}(b)$ & $1,0$ &2&&&&&&&&&\\
 & $0,0$ &&&&&&&&&& \\
\hline
&\multirow{2}{*}{}$1,1$ &2&2&&&&&&&& \\
 & $0,1$ &&2&2&&&&&&& \\
 $N_4^{\al, \be}(b)$& $1,0$ &2&&&&&&&&& \\
 & $0,0$ &&2&&&&&&&& \\
\hline
&\multirow{2}{*}{}$1,1$ &2&6&&&&&&&& \\
 & $0,1$ &&4&2&2&&&&&& \\
 $N_5^{\al, \be}(b)$& $1,0$ &2&6&&&&&&&&$2 \leq N_{5}^{0,0}(4) \leq 2$\\
 & $0,0$ &&2&2&&&&&&& \\
\hline
&\multirow{2}{*}{}$1,1$ &2&12&2&&&&&&&  \\
 & $0,1$ &&8&4&2&2&&&&& \\
 $N_6^{\al, \be}(b)$& $1,0$ &2&10&&&&&&&& $2 \leq N_{6}^{0,0}(4) \leq 4$\\
 & $0,0$ &&4&6&2&&&&&& \\
\hline
&\multirow{2}{*}{}$1,1$ &2&24&6&&&&&&& \\
 & $0,1$ &&14&10&4&2&2&&&& \\
 $N_7^{\al, \be}(b)$& $1,0$ &2&16&14&&&&&&& $4 \leq N_{7}^{0,0}(4) \leq 8$\\
 & $0,0$ &&8&8&6&2&&&&& \\
\hline
&\multirow{2}{*}{}$1,1$ &2&46&14&2&&&&&&  \\
 & $0,1$ &&26&22&8&4&2&2&&& \\
 $N_8^{\al, \be}(b)$& $1,0$ &2&36&18&&&&&&&$6 \leq N_{8}^{0,0}(4) \leq 12$\\
 & $0,0$ &&12&24&10&6&2&&&& \\
\hline
&\multirow{2}{*}{}$1,1$ &2&86&34&6&&&&&&  \\
 & $0,1$ &&48&46&18&8&4&2&2&& \\
 $N_9^{\al, \be}(b)$& $1,0$ &2&66&38&22&&&&&&$8 \leq N_{9}^{0,0}(4) \leq 26$\\
 & $0,0$ &&22&50&20&10&6&4&&& \\
\hline
&\multirow{2}{*}{}$1,1$ &2&160&78&14&2&&&&&  \\
 & $0,1$ &&88&96&40&16&8&4&2&2&\\
 $N_{10}^{\al, \be}(b)$& $1,0$ &2&110&110&18&&&&&&$12 \leq N_{10}^{0,0}(4) \leq 46$ \\
 & $0,0$ &&42&94&60&24&12&6&2&&\\
%\hline
%&\multirow{2}{*}{}$1,1$ &2&296&174&34&6&&&&&&  \\
% & $0,1$ &&164&202&80&34&16&8&4&2&2& \\
 %$N_{11}^{\al, \be}(b)$& $1,0$ &2&206&210&64&30&&&&&&$18 \leq N_{11}^{0,0}(4) \leq 86$ \\
 %& $0,0$ &&72&206&110&48&24&12&6&2&& \\
\hline
%\end{tabular}
\end {tabular}
\end{center}

\end {table}

\newpage

\bibliographystyle{plain}
\bibliography{m1}

\begin{thebibliography}{10}

\bibitem{Bru:Ross}
Richard~A. Brualdi and Jeffrey~A. Ross.
\newblock On the exponent of a primitive, nearly reducible matrix.
\newblock {\em Math. Oper. Res.}, {\bf 5}({\bf 2}):229--241, 1980.

\bibitem{bru}
Richard~A. Brualdi and Herbert~J. Ryser.
\newblock {\em Combinatorial matrix theory}, volume~{\bf 39} of {\em
  Encyclopedia of Mathematics and its Applications}.
\newblock Cambridge University Press, New York, 2013.
\newblock Paperback edition of the 1991 original [ MR1130611].

\bibitem{D:S}
Charles Delorme and Patrick Sol\'e.
\newblock Diameter, covering index, covering radius and eigenvalues.
\newblock {\em European J. Combin.}, {\bf 12}({\bf 2}):95--108, 1991.

\bibitem{D:M}
A.~L. Dulmage and N.~S. Mendelsohn.
\newblock Gaps in the exponent set of primitive matrices.
\newblock {\em Illinois J. Math.}, {\bf8}:642--656, 1964.

\bibitem{Ho:Va}
John~C. Holladay and Richard~S. Varga.
\newblock On powers of non-negative matrices.
\newblock {\em Proc. Amer. Math. Soc.}, {\bf 9}:631--634, 1958.

\bibitem{Kim:Song:Hwang}
Byeong~Moon Kim, Byung~Chul Song, and Woonjae Hwang.
\newblock Nonnegative primitive matrices with exponent 2.
\newblock {\em Linear Algebra Appl.}, {\bf 407}:162--168, 2005.

\bibitem{Lewin}
Mordechai Lewin.
\newblock Bounds for exponents of doubly stochastic primitive matrices.
\newblock {\em Math. Z.}, {\bf 137}:21--30, 1974.

\bibitem{M:Y}
Mordechai Lewin and Yehoshua Vitek.
\newblock A system of gaps in the exponent set of primitive matrices.
\newblock {\em Illinois J. Math.}, {\bf 25}({\bf 1}):87--98, 1981.

\bibitem{Lic:Cai}
Huangfu Lichao and Junliang Cai.
\newblock The characterization of symmetric primitive matrices with exponent
  {$n-3$}.
\newblock {\em Int. J. Math. Comb.}, {\bf 3}:56--78, 2008.

\bibitem{Liu1}
Bo~Lian Liu.
\newblock A note on the exponents of primitive {$(0,1)$} matrices.
\newblock {\em Linear Algebra Appl.}, {\bf 140}:45--51, 1990.

\bibitem{B:M:W}
Bo~Lian Liu, Brendan~D. McKay, Nicholas~C. Wormald, and Ke~Min Zhang.
\newblock The exponent set of symmetric primitive {$(0,1)$} matrices with zero
  trace.
\newblock {\em Linear Algebra Appl.}, {\bf 133}:121--131, 1990.

\bibitem{Liuu:Youu:Yuu}
Bolian Liu, Lihua You, and Gexin Yu.
\newblock On extremal matrices of second largest exponent by {B}oolean rank.
\newblock {\em Linear Algebra Appl.}, {\bf 422}({\bf 1}):186--197, 2007.

\bibitem{minc}
Henryk Minc.
\newblock {\em Nonnegative matrices}.
\newblock Wiley-Interscience Series in Discrete Mathematics and Optimization.
  John Wiley \& Sons, Inc., New York, 1988.
\newblock A Wiley-Interscience Publication.

\bibitem{M:S}
Monimala Nej and Satyanarayana~A. Reddy.
\newblock Binary strings of length $n$ with $x$ zeros and longest $k$-runs of
  zeros.
\newblock {\em arXiv:1707.02187}, 2017.

\bibitem{NY1}
M.~A. Nyblom.
\newblock Enumerating binary strings without {$r$}-runs of ones.
\newblock {\em Int. Math. Forum}, {\bf 7}({\bf 37-40}):1865--1876, 2012.

\bibitem{Per}
Peter Perkins.
\newblock A theorem on regular matrices.
\newblock {\em Pacific J. Math.}, {\bf 11}:1529--1533, 1961.

\bibitem{Ross}
Jeffrey~A. Ross.
\newblock On the exponent of a primitive, nearly reducible matrix. {II}.
\newblock {\em SIAM J. Algebraic Discrete Methods}, {\bf 3}({\bf 3}):395--410,
  1982.

\bibitem{HS}
Hans Schneider.
\newblock Wielandt's proof of the exponent inequality for primitive nonnegative
  matrices.
\newblock {\em Linear Algebra Appl.}, {\bf 353}:5--10, 2002.

\bibitem{shao}
Jia~Yu Shao.
\newblock The exponent set of symmetric primitive matrices.
\newblock {\em Sci. Sinica Ser. A}, {\bf 30}({\bf 4}):348--358, 1987.

\bibitem{F:M:J}
Fuyi Wei, Maoquan Chu, and Jianzhong Wang.
\newblock On odd primitive graphs.
\newblock {\em Australas. J. Combin.}, {\bf19}:11--15, 1999.

\bibitem{Wie}
Helmut Wielandt.
\newblock Unzerlegbare, nicht negative {M}atrizen.
\newblock {\em Math. Z.}, {\bf 52}:642--648, 1950.

\bibitem{ZKM}
Ke~Min Zhang.
\newblock On {L}ewin and {V}itek's conjecture about the exponent set of
  primitive matrices.
\newblock {\em Linear Algebra Appl.}, {\bf 96}:101--108, 1987.

\end{thebibliography}

\end{document}